\documentclass{amsart}
\title[Linear transformations preserving the $PF$-property]
{On linear transformations preserving the P\'olya frequency property}

\author{Petter Br\"and\'en}

\thanks{
Research financed by the EC's IHRP Programme, within the Research Training 
Network 
``Algebraic Combinatorics in Europe'', grant HPRN-CT-2001-00272, while the 
author was at Universit\'a di Roma ``Tor Vergata'', Rome, Italy.} 

\address{Matematik,
  Chalmers tekniska h\"ogskola och G\"oteborgs universitet,\linebreak
  S-412~96  G\"oteborg, Sweden}
\email{branden@math.chalmers.se}
\date{March 22, 2004}

\usepackage[english]{babel}
\usepackage{amsmath}
\usepackage{amsthm,amssymb,latexsym}
\usepackage{t1enc}
\usepackage{epic}
\usepackage{eepic}
\usepackage{mathrsfs}

\numberwithin{equation}{section}
\newtheorem{proposition}{Proposition}[section]
\newtheorem{lemma}[proposition]{Lemma}
\newtheorem{corollary}[proposition]{Corollary}
\newtheorem{theorem}[proposition]{Theorem}
\newtheorem{conjecture}[proposition]{Conjecture}
\theoremstyle{definition}
\newtheorem{definition}[proposition]{Definition}

\newcommand{\N}{\mathbb{N}}

\newcommand{\Z}{\mathbb{Z}}

\newcommand{\Four}{\mathcal{L}_{\phi}}
\newcommand{\E}{\mathcal{E}}
\newcommand{\Ref}{\mathcal{R}}
\newcommand{\R}{\mathbb{R}}

\newcommand{\sym}{\mathcal{S}}

\newcommand{\Al}{\mathscr{A}}

\def\hy#1,#2,#3,#4,{{}_2F_1\left(#1,#2;#3;#4\right)}
\def\emm#1,{{\em #1}}
\def\ba#1,{\overline{#1}}

\def\gen#1,{\langle #1 \rangle}
\def\geno#1,{\langle #1 \rangle_{\infty}}

\def\newop#1{\expandafter\def\csname #1\endcsname{\mathop{\rm
#1}\nolimits}}

\newop{span}
\newop{min}
\newop{max}
\newop{deg}
\newop{des}
\newop{exc}
\newop{c}
\newop{Con}
\newop{Ann}
\newop{dim}
\newop{Int}
\newop{Im}
\newop{ker}
\newop{mult}
\newop{rank}
\newop{sign}
\newop{sgn}
\newop{TrCl}
\newop{inf}

\usepackage{graphicx}

\begin{document}
\maketitle

\begin{center}
{\em This is the final version to appear in Trans. Amer. Math. Soc.}
\end{center}

\begin{abstract}
We prove that certain linear operators preserve the P\'olya frequency 
property and real-rootedness, and apply our results to settle some 
conjectures and open problems in 
combinatorics proposed by B\'ona, Brenti and Reiner-Welker.
\end{abstract}
\nocite{*}\bibliographystyle{plain}
\thispagestyle{empty}
\section{Introduction} 
Many sequences encountered in various areas of mathematics, statistics and computer 
science are known or conjectured to be unimodal or log-concave, see 
\cite{brentisurvey,stanleysurvey,stanleypositivity}. A sufficient 
condition for a sequence to enjoy these properties is that it is 
a P\'olya frequency ($PF$ for short) sequence, or equivalently for finite sequences, that its 
generating function has only real and non-positive zeros. It is often the 
case that the generating function of a finite $PF$-sequence has more 
transparent properties when expanded in a basis other than the standard basis 
$\{x^i\}_{i \geq 0}$ of $\R[x]$. Therefore it is natural to investigate how 
$PF$-sequences translate when expressed in various basis. This amounts 
to studying properties of the linear operator that maps one basis to another. 
A systematic study of this was first pursued by Brenti in \cite{brentithesis}. 
This is also the theme of this paper. 
 
In Section \ref{theory} we will study linear operators of the type 
$$
\phi_F = \sum_{k=0}^n Q_k(x) \frac {d^k} {dx^k},
$$     
where $F(x,z) = \sum_{k=0}^n Q_k(x)z^k \in \R[x,z]$.  
Here we will give 
sufficient conditions on $F$ for $\phi_F$ to preserve the $PF$-property. The results attained 
generalizes and unifies theorems of Hermite, Poulain, P\'olya and Schur.
We will also in this section give a sufficient condition for a family of natural $\R$-bilinear 
forms to preserve the $PF$-property in both arguments. This generalizes results 
of Wagner \cite{brentietal,wagnerpp,wagnernc}. 

An important linear operator in combinatorics is the operator defined 
by $\E( \binom x i ) = x^i$, for all $i \in \N$. In Section \ref{special}  we will prove 
that whenever a polynomial $f$ of degree $d$ has nonnegative coefficients when expanded in 
the basis $\{ x^i(x+1)^{d-i} \}_{i=0}^d$ the polynomial $\E(f)$ will have only 
real, non-positive and simple zeros.    

In the remainder of the paper we use the theory developed to settle some conjectures and 
open problems raised in combinatorics.  In Section \ref{tstack} we prove 
that the numbers $\{W_t(n,k)\}_{k=0}^{n-1}$ of $t$-stack sortable permutations 
in $\sym_n$ with $k$ descents form $PF$-sequences when $t = 2,n-2$, and thereby  
settling two new cases of an open problem proposed by B\'ona 
\cite{bonastack,bonasym}. 

In Section \ref{qeuler} we prove that the 
$q$-Eulerian polynomials, $A_n(x;q)$, defined by Foata and Sch\"utzenberger \cite{foata} 
and further studied by Brenti in \cite{brentiplethysm} have only real zeros for all 
integers $q$. This 
settles a conjecture 
raised by Brenti. Here we also continue the study of the $W$-Eulerian polynomials, defined 
for any finite Coxeter group $W$ and the $q$-analog $B_n(x;q)$, initiated by 
Brenti in \cite{brentieuler}. 

In Section \ref{hvector} we prove that the $h$-vectors of a family simplicial complexes 
associated to finite Weyl groups defined by Fomin and Zelevinski \cite{fomin} are $PF$, thus 
settling an open problem raised by Reiner and Welker \cite{reinerwelker}.  

\section{Notation and preliminaries}\label{notation}
In this section we collect definitions, notation and results that will be used 
frequently in the rest of the paper.  
Let $\{a_i\}_{i=0}^\infty$ be a sequence of real numbers. It is {\em unimodal} if 
there is a number $p$ such that $a_0 \leq a_1 \leq \cdots \leq a_p \geq a_{p+1} \geq \cdots$, 
and  {\em log-concave} if  
$a_i^2 \geq a_{i-1}a_{i+1}$ for all $i > 0$. 

An infinite matrix $A= (a_{ij})_{i,j \geq 0}$ of real numbers is {\em totally positive}, 
$TP$, if all minors of $A$ are nonnegative. An infinite sequence 
$\{a_i\}_{i=0}^\infty$ of real numbers is a {\em P\'olya frequency sequence}, 
$PF$-sequence,  if 
the matrix $(a_{i-j})_{i,j \geq 0}$ is $TP$. Thus a $PF$-sequence is by definition 
log-concave and therefore also unimodal. A finite 
sequence $a_0, a_1, a_2, \ldots, a_n$ is said to be $PF$ if the infinite 
sequence $a_0, a_1, a_2, \ldots, a_n,0,0, \ldots$ is $PF$. A sequence $\{a_i\}_{i=0}^{\infty}$ 
is said to be $PF_r$ if all minors of size at most $r$ of  $(a_{i-j})_{i,j \geq 0}$ are 
nonnegative. If the polynomials 
$\{b_i(x)\}_{i=0}^d$ are linearly independent over $\R$ and $r \in \N$ we define the 
set $PF_r[\{b_i(x)\}_{i=0}^d]$ to be 
$$
PF_r[\{b_i(x)\}_{i=0}^d]= \{ \sum_{i=0}^d \lambda_ib_i(x) : \{\lambda_i\}_{i=0}^{\infty} 
\mbox{ is } PF_r\}, 
$$
and $PF[\{b_i(x)\}_{i=0}^d]= \bigcap_{r=1}^{\infty}PF_r[\{b_i(x)\}_{i=0}^d]$.
 
The following theorem characterizes $PF$-sequences. It was conjectured by Schoenberg and 
proved by Edrei \cite{edrei}, see also \cite{karlin}. 
\begin{theorem}
Let $\{a_i\}_{i=0}^{\infty}$ be a sequence of real numbers with $a_0=1$. Then 
it is a $PF$-sequence if and only if the generating function 
can be expanded, in a neighborhood of the origin, as 
$$
\sum_{i \geq 0}a_iz^i = e^{\gamma z} \frac { \prod_{i\geq 0}(1+\alpha_iz)}
                                           { \prod_{i\geq 0}(1-\beta_iz)},
$$
where $\gamma \geq 0$, $\alpha_i, \beta_i > 0$ and 
$\sum_{i\geq 0}(\alpha_i + \beta_i) < \infty$.
\end{theorem}
A consequence of this theorem is that a finite sequence is $PF$ if and only if 
its generating function is a polynomial with only real non positive zeros. 

Let $f, g \in \R[x]$ be real-rooted with zeros:
$\alpha_1 \leq \cdots \leq \alpha_i$ and 
$\beta_1 \leq \cdots \leq \beta_j$, respectively. We say that $f$ {\em interlaces} $g$, denoted 
$f \preceq g$, if $j=i+1$ and 
$$
\beta_1 \leq \alpha_1 \leq \beta_2 \leq \cdots \leq \beta_{j-1} \leq \alpha_{j-1} \leq \beta_j. 
$$
We say that $f$ {\em alternates left of} $g$, denoted $f \ll g$,  if $i=j$ and 
 $$
\alpha_1 \leq \beta_1 \leq \cdots \leq \beta_{i-1} \leq \alpha_{i} \leq \beta_i. 
$$
If in addition $f$ and $g$ have no common zero then we say that $f$ 
{\em strictly interlaces} $g$ and 
$f$ {\em strictly alternates left of} $g$, respectively. 
 We also say that two polynomials $f$ and $g$ {\em alternate} if one of 
the polynomials alternates left of or interlaces the other. 
We will need two simple lemmata 
concerning these concepts. A polynomial is said to be {\em standard} if its leading 
coefficient is positive.  
\begin{lemma}\label{positivesum}
Let $g$ and $\{f_i\}_{i=1}^n$ be real-rooted standard polynomials. 
\begin{itemize}
\item[(i)] If for each $1 \leq i \leq n$ 
we have either $g \ll f_i$ or $g \preceq f_i$. Then the sum 
$F=f_1 + f_2 + \cdots + f_n$ is real-rooted with $g \preceq F$ or $g \ll F$, depending 
on the degree of $F$. 
\item[(ii)] If for each $1 \leq i \leq n$ 
we have either $f_i \ll g$ or $f_i \preceq g$. Then the sum 
$F=f_1 + f_2 + \cdots + f_n$ is real-rooted with $F \preceq g$ or $F \ll g$, depending 
on the degree of $F$. 
\end{itemize}

\end{lemma} 
\begin{proof}
The lemma follows easily by counting the sign-changes of $F$ at the zeros of $g$, 
see e.g., \cite[Prop. 3.5]{wagnertot}. 
\end{proof}
The next lemma is obvious:
\begin{lemma}\label{weaktransitivity}
If $f_0, f_1, \ldots, f_n$ are real-rooted polynomials with $f_0 \ll f_n$ and 
$f_{i-1} \ll f_i$ for all $1 \leq i \leq n$, then $f_i \ll f_j$ for 
all $0 \leq i \leq j \leq n$. 
\end{lemma}
The following theorem is a characterization of alternating polynomials due to 
Obreschkoff \cite{obreschkoff} and Dedieu \cite{dedieu}:
\begin{theorem}
Let $f, g \in \R[x]$. Then $f$ and $g$ alternate (strictly alternate) if 
and only if all polynomials in the space
$$
\{\alpha f + \beta g : \alpha,\beta \in \R\},
$$
have only real (real and simple) zeros.
\end{theorem}
An immediate but non-trivial consequence of this theorem is:
\begin{corollary}\label{obreschkor}
Let $\phi : \R[x] \rightarrow \R[x]$ be a linear operator. Then 
$\phi$ preserves the real-rootedness property (real- and simple-rootedness property) 
only if $\phi$ preserves the alternating property 
(strictly alternating property).
\end{corollary} 
  
We denote by $\N$ the set of natural numbers 
$\{0,1,2,\ldots\}$. The symmetric group of bijections 
$\pi : \{1,2, \ldots, n\} \rightarrow \{1,2, \ldots, n\}$ is denoted by 
$\sym_n$. A {\em descent} in a permutation $\pi \in \sym_n$ is an index 
$1 \leq i \leq n-1$ such 
that $\pi(i) > \pi(i+1)$. Let $\des(\pi)$ denote the number of descents in $\pi$. The 
{\em Eulerian polynomials}, $A_n(x)$,  are defined by 
$A_n(x)= \sum_{\pi \in \sym_n}x^{\des(\pi)+1}$ and satisfies, see e.g. \cite{comtet} 
$$
\sum_{k\geq 0} k^n x^k = \frac {A_n(x)}{(1-x)^{n+1}}.
$$   
The binomial polynomials are defined by $\binom x 0 = 1$ and 
$\binom x k = \frac {x(x-1)\cdots (x-k+1)}{k!}$ for $k \geq 1$.

In several proofs we will implicitly use the fact that the zeros of a polynomial are 
continuous functions of the coefficients of the polynomial. In particular 
the limit of a sequence of real-rooted polynomials is again real-rooted. For 
a treatment of these matters we refer the reader to \cite{marden}.

\section{A class of linear operators preserving the $PF$-property}\label{theory}
For any polynomial $F(x,z) = \sum_{k=0}^nQ_k(x)z^k \in \R[x,z]$ we define 
a linear operator $\phi_F : \R[x] \rightarrow \R[x]$ by, 
$$
\phi_F(f) := \sum_{k=0}^nQ_k(x)\frac {d^k} {dx^k}f(x).
$$
In this section we will investigate for which $F \in \R[x,z]$ the linear 
operator $\phi_F$ preserves real-rootedness- and the $PF$-property .   

We will need some terminology and a theorem from \cite{branden1}. 
For $\xi \in \R$ let $T_{\xi} : \R[x] \rightarrow \R[x]$ be the 
translation operator defined by $T_{\xi}(f(x)) = f(x+\xi)$.  
For any linear operator $\phi : \R[x] \rightarrow \R[x]$ we define a 
linear transform $\Four : \R[x] \rightarrow \R[x,z]$ by 
\begin{eqnarray}\label{defour}
\Four(f) &:=& \phi(T_z(f)) \nonumber \\
         &=& \sum_n \phi(f^{(n)})(x) \frac {z^n} {n!} \\
         &=& \sum_n \frac {\phi(x^n)}{n!} f^{(n)}(z). \nonumber
\end{eqnarray}  

\begin{definition}\label{al}
Let $\phi : \R[x] \rightarrow \R[x]$ be a linear operator. 
Define a function $d_\phi: \R[x] \rightarrow \N \cup \{-\infty\}$ by: 
If $\phi(f^{(n)})=0$ for all 
$n \in \N$, we let $d_\phi(f)=-\infty$. Otherwise let $d_\phi(f)$ be the 
smallest integer $d$ such that $\phi(f^{(n)})=0$ for all $n > d$. 
Hence $d_{\phi}(f) \leq \deg f$ for all $f \in \R[x]$.

The set $\Al(\phi)$ is defined as follows: If 
$d_\phi(f) \in \{-\infty,0\}$ and $\phi(f)$ is standard real- and simple-rooted, then 
$f \in \Al(\phi)$. Moreover, $f \in \Al(\phi)$ if  
$d=d_\phi(f) \geq 1$ and all of the following conditions are satisfied: 
\begin{itemize}
\item[(i)] $\phi(f^{(i)})$ all have leading coefficients of the same sign and 
$\deg(\phi(f^{(i-1)}))= \deg(\phi(f^{(i)}))+1$ 
for $1 \leq i \leq d$,
\item[(ii)] $\phi(f)$ and $\phi(f')$ have no common real zero,
\item[(iii)] $\phi(f^{(d)})$ strictly interlaces $\phi(f^{(d-1)})$,
\item[(iv)] for all $\xi \in \R$ the polynomial 
$\Four(f)(\xi,z)$ is real-rooted.
\end{itemize}
\end{definition} 
The following theorem is proved in \cite{branden1}:
\begin{theorem}\label{huvud}
Let $\phi : \R[x] \rightarrow \R[x]$ be a linear operator. If 
$f \in \Al(\phi)$ then $\phi(f)$ is real- and simple-rooted. 

\end{theorem}
We will also need the following classical theorem of Hermite and Poulain. For a proof 
see \cite{obreschkoff}. 
\begin{theorem}\label{hp}
Let $f=a_0 +a_1x+ \cdots + a_nx^n$ and $g$ be real-rooted polynomials. 
Then the polynomial 
$$
f(\frac d {dx})g:= a_0g(x) + a_1g'(x) + \cdots + a_ng^{(n)}(x)
$$
is real-rooted. Moreover, if $f(\frac d {dx})g \neq 0$ then any multiple zero 
of $f(\frac d {dx})g$ is a multiple zero of $g$. 
\end{theorem}
 The following theorem gives a sufficient condition for a polynomial to be mapped 
onto a real-rooted polynomial. 
\begin{theorem}\label{hyper}
Let $F=\sum_{k=0}^nQ_k(x)z^k$ be such that $Q_0 \neq 0$ and  
\begin{itemize}
\item[(I)] For all $\xi \in \R$, $F(\xi,z)$ is real-rooted, 
\item[(II)] $Q_0$ strictly interlaces or strictly alternates left of $Q_1$, 
and $\deg Q_0 =0$ or $Q_0$ and $Q_1$ have leading coefficients of the same sign. 
\end{itemize}
Suppose that 
\begin{itemize}
\item[ (III)] $f$ is real- and simple-rooted and 
that for $0 \leq k \leq \deg f$ the 
polynomials $\phi_F(f^{(k)})$ have their leading term of the same sign with 
$$
\deg \phi_F(f^{(k)}) = \deg Q_0 + \deg f - k.
$$
\end{itemize}
Then $\phi_F(f)$ is real- and simple-rooted.
\end{theorem}

\begin{proof}
We will show that the set of real- and simple-rooted polynomials satisfying 
(III) is a subset of $\Al(\phi_F)$ by verifying conditions (i)-(iv) of 
Definition \ref{al}. Condition (i) follows immediately from (III). For 
condition (iv) note that 
$$
\Four(f)(\xi,z) = \sum_{k=0}^nQ_k(\xi)f^{(k)}(\xi+z),
$$
so by the Theorem \ref{hp} condition (iv) is satisfied. 
Suppose that $\eta$ is a common zero of $\phi_F(f)$ and $\phi_F(f')$. By 
\eqref{defour} we have that $0$ is a multiple zero of 
$\Four(f)(\eta,z)$. Moreover, since  $\Four(f)(\eta,z)$ is not identically equal to zero, by 
(II), Theorem \ref{hp} 
tells us that $0$ is a multiple zero of $f(\eta +z)$. This means that $\eta$ is a multiple 
zero of $f$ contrary to the assumption that $f$ is simple-rooted, and verifies condition (ii). 

For condition (iii) we have to show that for all $\alpha \in \R$ such that 
$x+\alpha$ satisfies (III) the 
polynomial $\phi_F(1)=Q_0$ strictly interlaces 
$f(x):=\phi_F(x+\alpha)=(x+\alpha)Q_0+Q_1$. This 
follows from (II) when analyzing the sign of $f(x):=\phi_F(x+\alpha)$ at the 
zeros of $Q_0$: Let $\alpha_k < \alpha_{k-1} < \cdots < \alpha_1$ be the 
zeros of $Q_0$ ordered by size. Suppose that $Q_0$ and $Q_1$ are standard and that $Q_0$ 
strictly interlaces or strictly alternates 
left of $Q_1$. Then $\sgn f(\alpha_i)=\sgn Q_1(\alpha_i)=(-1)^i$ for $1 \leq i \leq k$. By 
Rolle's theorem we know that $f$ has a zero in each interval 
$(\alpha_i, \alpha_{i+1})$. This accounts for $k-1$ real zeros of $f$. Since 
$Q_0$ has positive sign, so does $f$ by condition (III). Now, 
because $f(\alpha_1)<0$ and $f$ is standard, $f$ must have a zero to 
the right of $\alpha_1$. We now know that $f$ has $k$ zeros 
real. The signs at $\alpha_i$ forces the remaining zero to be in the interval 
$(-\infty, \alpha_k)$. Thus $Q_0$ strictly interlaces $f$ as was to be 
shown. 

Now, if $Q_0=A \in \R$ then $\deg Q_1  \leq 1$. Suppose that $Q_1=B \in \R$. 
Then clearly $A$ strictly interlaces $(x+\alpha)A +B$. If 
$Q_0=A$ and $Q_1=Cx+D$ where $A,C,D \in \R$, then $f=(A+C)x + A\alpha+D$, so by 
(III) we have that $Q_0$ strictly interlaces $f$. This concludes the proof.    
\end{proof}
In some cases it may be convenient to have sharper hypothesis. Therefore 
we state the following form of the theorem. 
\begin{corollary}\label{maincor}
Let $d \in \N$ be given and let $F=\sum_{k=0}^nQ_k(x)z^k$ be such that 
$Q_0 \neq 0$ and  
\begin{itemize}
\item[(i)] For all $\xi \in \R$, $F(\xi,z)$ is real-rooted, 
\item[(ii)] $Q_0$ strictly interlaces or strictly alternates left of $Q_1$, 
and $\deg Q_0 =0$ or $Q_0$ and $Q_1$ have leading coefficients of the same sign. 
\item[ (iii)]  The polynomials $\phi_F(x^k)$, $0 \leq k \leq d$ have the same 
sign and     
$$
\deg \phi_F(x^k) = \deg Q_0 + k.
$$
\end{itemize}
Then $\phi_F(f)$ is real-rooted (real- and simple-rooted) if 
$f$ is real-rooted (real and simple-rooted) and $\deg(f) \leq d$. 
\end{corollary}
\begin{proof}
The case of real- and simple-rooted $f$ follows immediately from Theorem 
\ref{hyper} since (iii) implies (III). If $f$ is a real-rooted polynomial 
of degree at most $d$, then $f$ is the limit of a sequence 
$\{f_k\}_{k=0}^{\infty}$ of real- and simple-rooted polynomials of degree 
at most $d$. It follows that $\phi_F(f)$ is the limit of $\phi_F(f_k)$, and 
the thesis follows by continuity.
\end{proof}
In the language of $PF$-sequences we have:
\begin{theorem}
Let $d \in \N$ be given and let $F= \sum_{k=0}^nQ_k(x)z^k \in \R[x,z]$ be such that 
$Q_0 \neq 0$ and   
\begin{itemize}
\item[(i)] For all $\xi \in \R$, $F(\xi,z)$ is real-rooted, 
\item[(ii)] $\phi_F(1)$ strictly interlaces $\phi_F(x)$. 
\item[ (iii)] For all $0 \leq k \leq d$    
$$
\deg \phi_F(x^k) = \deg Q_0 + k,
$$
and $\phi_F(x^k) \in PF_1$.
\end{itemize}
Then $PF[\{\phi_F(x^i)\}_{i=0}^d] \subseteq PF[x^i]$.
\end{theorem}
Several old results can be derived from these last few theorems. In 
\cite[p. 163]{polyagen} P\'olya gave a theorem which he states probably was the most general 
theorem on real-rootedness known at the time. ''Dieser Satz geh\"ort wohl zu den 
allgemeinsten bekannten S\"atzen \"uber Wurzelrealit\"at.'':
\begin{theorem}
Let $f(x)$ be a real-rooted polynomial of degree $n$, and let 
$$
b_0 + b_1x + \cdots + b_{n+m}x^{n+m}, \ \ (m \geq 0)
$$
be a real-rooted polynomial such that $b_i>0$ for $0 \leq i \leq n$. Then 
the equation
$$
G(x,y) := b_0f(y) + b_1xf'(y) + b_2x^2f^{''}(y) + \cdots + b_nx^nf^{(n)}(y)=0,
$$
has $n$ real intersection points, (counted with multiplicity), with 
the line
$$
sx-ty+u =0,
$$
provided that $s, t \geq 0$, $s+t >0$ and $u \in \R$. 
\end{theorem}
\begin{proof} 
We may assume that $s,t >0$ since the other cases follows by continuity 
when $s$ and/or $t$ tends to zero. Thus we may write the equation as 
$$
a_0g(x) + a_1xg'(x) + a_2x^2g''(x) + \cdots + a_nx^ng^{(n)}(x)=0,
$$
where $g(x) = f(st^{-1}x+ut^{-1})$ and $a_i=s^it^{-i}b_i$. Now, we see   
that all hypothesis of Corollary \ref{maincor} are satisfied for 
$$
F(x,z)= a_0 + a_1xz + a_2x^2z^2 + \cdots + a_{n+m}x^{n+m}z^{n+m},
$$ 
when $d =n$. 
\end{proof}
We will later need one famous consequence of this theorem, $t=1,s=u=0$, due to Schur 
\cite{schur}. 
\begin{theorem}\label{schurmult}
Let $f=\sum_{k=0}^{n}a_kx^k$ and $g=\sum_{k=0}^{m}b_kx^k$ be two 
real-rooted polynomials such that $g$ has all zeros of the same sign. Then 
the polynomial 
$$
(f S g)(x) = \sum_{k \geq 0}^M k!a_kb_kx^k,
$$
where $M=\min(m,n)$ has only real zeros. 
\end{theorem}
\subsection{Multiplier-sequences}
A multiplier-sequence is a sequence $T=\{\gamma_i\}_{i=0}^{\infty}$ of real numbers 
such that if a polynomial $f(x)=a_0 +a_1x + \cdots + a_nx^n$ has only 
real zeros, then the polynomial
$$
T[f(x)]:= a_0\gamma_0 +a_1\gamma_1x + \cdots + a_n\gamma_nx^n,
$$
also has only real zeros. There is a characterization of multiplier-sequences 
due to P\'olya and Schur \cite[p. 100-124]{polyagen}:
\begin{theorem}
Let $T=\{\gamma_i\}_{i=0}^{\infty}$ be a  sequence of non-negative real 
numbers and let 
$\phi(x) = T[e^x] = \sum_{k=0}^{\infty}\gamma_k\frac{x^k}{k!}$ be its 
exponential 
generating function. Then $T$ is 
a multiplier-sequence if and only if $\phi$ is a real entire function which 
can be written as  
$$
\phi(x)=cx^ne^{\beta x} \prod_{k=1}^{\infty}(1+\delta_k x),
$$
where $c>0$, $\beta \geq 0$, $\delta_k \geq 0$, 
$n \in \N$ and $\sum_{k=1}^{\infty}\delta_k < \infty$. 
\end{theorem}

The following lemma is well-known but elementary, so we give a proof here. 
\begin{lemma}\label{internal}
A multiplier-sequence is strictly log-concave.  
In particular, a nonnegative multiplier-sequence has no internal zeros.
\end{lemma}
\begin{proof}
If $f(x)= a_mx^m + a_{m+1}x^{m+1} + \cdots + a_nx^n$ is real-rooted with $a_ma_n \neq 0$, 
then the coefficients satisfy (see \cite[p. 52]{hardy}):
$$
\frac {a_i^2} {{\binom n i}^2} > \frac {a_{i-1}}{\binom n {i-1}}\frac {a_{i+1}}{\binom n {i+1}}
\ \ (m < i < n).
$$
Now, if $\Gamma = \{\gamma_i\}_{i=0}^{\infty}$ is a multiplier-sequence, then 
$\Gamma[(x+1)^n]$ is real-rooted for all $n \in \N$, which implies 
$$
\gamma_i^2 > \gamma_{i-1} \gamma_{i+1},
$$
for all $i$ such that there are integers $m < i < n$ with $\gamma_m \gamma_n \neq 0$.  

\end{proof}  
\begin{theorem}\label{products}
Let $\{\lambda_k\}_{k=0}^{\infty}$ be a non-negative multiplier-sequence, 
and let $\alpha < \beta \in \R$ be given.  
Define two $\R$-bilinear forms $\R[x] \times \R[x] \rightarrow \R[x]$ by 
\begin{eqnarray*}
f\cdot g &:=& 
\sum_{k \geq 0} \frac {\lambda_k}{k!} f^{(k)}(x)g^{(k)}(x)(x-\alpha)^k(x-\beta)^k,\\
f \circ g &:=& \sum_{k \geq 0} \frac {\lambda_k}{k!} f^{(k)}(x)g^{(k)}(x)(x-\alpha)^k, 
\end{eqnarray*}
If $f$ is real-rooted and $g$ is $[\alpha,\beta]$-rooted, then $f\cdot g$ is 
real-rooted. If $f$ is real-rooted and $g$ is $[-\infty,\alpha]$-rooted, 
then $f\circ g$ is real-rooted.
\end{theorem}

\begin{proof}
We prove the statement for $\cdot$ since the case of $\circ$ is similar. We may assume 
that $\lambda_0 >0$. Clearly the theorem is true if $\lambda_i=0$ for all $i>0$, so 
by Lemma \ref{internal} we may assume that $\lambda_1>0$.  Let 
$g$ have all zeros simple and in the interval $(\alpha,\beta)$, and let
$\phi$ be the linear operator defined by
$\phi(f) = f \cdot g$. Then $\phi = \phi_F$, where 
$$
F(x,z)= \sum_{k\geq 0}\lambda_k \frac{g^{(k)}(x)}{k!}(x-\alpha)^k(x-\beta)^k z^k. 
$$ 
Since $\{\lambda_k\}_{k\geq 0}$ is a multiplier sequence $F(\xi,z)$ is real-rooted 
for all real choices of $\xi$. Now, $Q_0 = \lambda_0g(x)$ and 
$Q_1= \lambda_1(x-\alpha)(x-\beta)g'(x)$, so $Q_0$ strictly interlaces $Q_1$. Moreover, 
$\deg \phi(x^k)= \deg Q_0 +k$ for all $k$, so all the hypothesis of Corollary 
\ref{maincor} are fulfilled. Since any $[\alpha,\beta]$-rooted polynomial 
is the limit of polynomials which are $(\alpha,\beta)$- and simple-rooted the 
thesis follows by continuity. 
\end{proof}

There are a few bilinear forms on polynomials that occur frequently in combinatorics.
Let $\# : \R[x] \times \R[x] \rightarrow \R[x]$ be defined by 
$$
(f \# g)(x) := \sum_{k \geq 0} f^{(k)}(x)g^{(k)}(x)\frac{x^k}{k!}.
$$ 
This product is important when analyzing how the the zeros of $\sigma$-polynomials behave 
under disjoint union of graphs, see \cite{brentietal}. 
\begin{theorem}
Let $f$ be real-rooted and let $g$ have only real zeros of the same sign. Then 
$f \# g$ is real-rooted.  
\end{theorem} 
\begin{proof}
The theorem follows from Theorem \ref{products}, since $\{1\}_{k=0}^{\infty}$ is trivially 
a multiplier-sequence.
\end{proof}
This generalizes a result of Wagner, who proved that $f \# g$ is real-rooted whenever 
$f$ and $g$ have only non-negative zeros, see \cite{brentietal,wagnerpp}.

The {\em diamond product} of two polynomials $f$ and $g$ is given by 
\begin{equation}\label{diamondid}
(f \diamond g)(x)= \sum_{k \geq 0} \frac{f^{(k)}(x)}{k!}\frac{g^{(k)}(x)}{k!}x^k(x+1)^k.
\end{equation}
This product is important in the theory of $(P,\omega)$-partitions and 
the Neggers-Stanley conjecture and was first studied by Wagner in \cite{wagnernc,wagnertot}, 
see also Section \ref{special} of this paper.
Applying Theorem \ref{products} with the multiplier-sequence $\{\frac{1}{k!}\}_{k\geq 0}$ we 
get:
\begin{theorem}
Let $f$ be real-rooted and let $g$ have all zeros in the interval $[-1,0]$. Then 
$f \diamond g$ is real-rooted.  
\end{theorem}
This was first proved by Wagner \cite{wagnertot} under the additional hypothesis that 
$f$ has all zeros in $[-1,0]$, and generalized by the present author in \cite{branden1}. 

A sequence of real numbers $\Gamma = \{\gamma_k\}_{k=0}^{\infty}$ is called 
a {\em multiplier $n$-sequence } if for any real-rooted polynomial 
$f= a_0 + a_1x + \cdots + a_nx^n$ of degree at most $n$ the polynomial 
$\Gamma[f] := a_0\gamma_0 + a_1\gamma_1x + \cdots + a_n\gamma_nx^n$ is 
real-rooted. There is a simple algebraic characterization of multiplier $n$-sequences 
\cite{mult}:
\begin{theorem}\label{multiplier}
Let $\Gamma=\{\gamma_k\}_{k=0}^{\infty}$ be a sequence of real numbers. Then 
$\Gamma$ is a multiplier $n$-sequence if and only if 
$\Gamma[(x+1)^n]$ is real-rooted with all its zeros of the same sign. 
\end{theorem}
 
Recall the definition of the 
{\em hypergeometric function} $_2F_1$: 
$$
\hy a,b,c,z, = \sum_{m=0}^\infty \frac{(a)_m(b)_mz^m}{(c)_m m!},
$$
where $(\alpha)_0=1$ and $(\alpha)_m=\alpha(\alpha+1)\cdots (\alpha+m-1)$ 
when $m \geq 1$. The Jacobi polynomial $P^{(\alpha,\beta)}_n(x)$ can 
be expressed as follows \cite[p. 254]{rainville}:
\begin{equation}\label{ulv1}
P^{(\alpha,\beta)}_n(x)=\frac {(1+\alpha)_n}{n!} \hy -n,1+\alpha+\beta+n,1+\alpha,\frac{1-x}2,,
\end{equation}
 
We need the following lemma: 
\begin{lemma}\label{factor}
Let $n$ be a positive integer and $r$ a non-negative real number. Then  
$\Gamma=\{\binom{-n-r}k\}_{k=0}^{\infty}$ is a multiplier $n$-sequence. 
\end{lemma}
\begin{proof}
Let $r > 0$. Then    
\begin{eqnarray*}
\Gamma[(x+1)^n] &=& \sum_{k=0}^n\binom {-n-r} k \binom n k x^k \\
                &=& \hy -n,n+r,1,x,   \\
                &=& P^{(0,r-1)}_n(1-2x),
\end{eqnarray*}
where the last equality follows from \eqref{ulv1}.
Since the Jacobi polynomials 
are known, see \cite{rainville}, to have all their zeros in 
$[-1,1]$ when $\alpha , \beta> -1$,  
we have that $\Gamma[(x+1)^n]$ has all its zeros in $[0,1]$. The case 
$r=0$ follows by continuity when we let $r$ tend to zero from above.  
\end{proof}
For any real number $q$ let $\Gamma_q := \{q+k\}_{k=0}^{\infty}$. The following 
Corollary was known already to Laguerre: 
\begin{corollary}\label{gammaq}
Let $n > 1$ be a positive integer. Then $\Gamma_q$ is a multiplier $n$-sequence if 
and only if $q \notin (-n,0)$. 
\end{corollary}
\begin{proof}
Let $q \in \R$ be given. 
We have to determine for which $n >1$ the zeros of 
$\Gamma_q[(x+1)^n]$ are all real and of the same sign. Now, 
$$
\Gamma_q[(x+1)^n]= (x+1)^{n-1}\{ (n+q)x + q \},
$$
and the theorem follows.
\end{proof}

\section{The E-transformation}\label{special}
The {\em E-transformation} is the invertible linear operator, 
$\E : \R[x] \rightarrow \R[x]$, defined by 
$$
\E( \binom x i ) =x^i,
$$
for all $i \in \N$. The $PF$-preserving properties of this linear operator was first studied 
in \cite{brentithesis} and later in \cite{wagnernc,wagnertot} and \cite{branden1}. 
It is important in the theory 
of $(P,\omega)$-partitions since it maps the order-polynomial of a labeled poset to the 
$E$-polynomial of the same labeled poset, see
\cite{brentithesis,wagnernc}. In, \cite{brentithesis} Brenti proved the following theorem. Let 
$\lambda(f)$ and $\Lambda(f)$ denote the smallest and the largest real zero of 
the polynomial $f$, respectively. 
\begin{theorem}\label{brentiomega}
Suppose that $f \in \R[x]$ has only real zeros and that 
$f(n)=0$ for all $n \in ([\lambda(f),-1]\cup [0,\Lambda(f)])\cap \Z$. Then 
$\E(f)$ has all zeros real and non-positive.
\end{theorem}
In this section we will prove the following theorem:
\begin{theorem}\label{strangepf}
For all $n \in \N$ we have 
$$
PF_1[\{x^i(x+1)^{n-i}\}_{i=0}^n] \subseteq PF[\binom x i] 
$$
Moreover if $f \in  PF_1[\{x^i(x+1)^{n-i}\}_{i=0}^n]$ then 
$\E(f)$ has simple zeros and 
$$
\E((x+1)^d) \ll \E(f) \ll \E(x^d).
$$
\end{theorem}
The diamond product \eqref{diamondid} is intimately connected with the E-transformation. 
By the Vandermonde identity 
$$
\binom x i \binom x j = \sum_{k \geq 0} \binom k {k-i,i+j-k,k-j} \binom x k,
$$
it follows, see \cite{wagnertot}, that 
\begin{equation}\label{diamonddef}
(f \diamond g)(x) := \E(\E^{-1}(f)\E^{-1}(g)).
\end{equation}
We will later need a symmetry property of $\E$. Let $\Ref : \R[x] \rightarrow \R[x]$ 
be the algebra automorphism defined by $\Ref(x) = -1-x$. 
\begin{lemma}\label{symE}
$$
\Ref \E = \E \Ref
$$
\end{lemma}
\begin{proof}
Let $n$ be a nonnegative integer. Using the identity 
$$
\binom{x+n}n = \sum_{k=0}^n \binom n k \binom x k,
$$
and the fact that $\binom {-x-1}n = (-1)^n\binom {x+n} n$
we get 
\begin{eqnarray*}
\E \Ref  \binom x n  &=&(-1)^n \E \binom {x+n} n \\
&=& (-1)^n \sum_{k=0}^n \binom n k \E\binom x k \\
&=& (-1-x)^n \\
&=& \Ref \E \binom x n,
\end{eqnarray*}
and the lemma follows.
\end{proof}

\begin{lemma}\label{rr}
Let $\alpha \in [-1,0]$ and let $f$ be a polynomial such that $\E(f)$ is $[-1,0]$-rooted. Then 
$\E((x-\alpha)f)$ is $[-1,0]$-rooted  and $\E(f)$ interlaces $\E((x-\alpha)f)$. 
If $\E(f)$ in addition only has 
simple zeros, then so does $\E((x-\alpha)f)$. 
\end{lemma}
\begin{proof}
Let $g=\E(f)$ and let $\alpha \in [-1,0]$. By \eqref{diamondid} and \eqref{diamonddef} 
we have that 
\begin{equation}\label{ylv}
\E((x-\alpha)f)= (x-\alpha)g + x(x+1)g'. 
\end{equation}
Since $g$ interlaces $(x-\alpha)g$ and $x(x+1)g'$ it also interlaces the sum, 
by Lemma \ref{positivesum}. Also, if $x \notin [-1,0]$ then 
the summands have the same sign so  $\E((x-\alpha)f)$ cannot have any zeros outside $[-1,0]$. 
Suppose that $g$ has only simple zeros. Then by \eqref{ylv} the only possible 
common zeros of $g$ and  $\E((x-\alpha)f)$ are $0$ and $-1$. If $\deg(f) \geq 1$ it 
also follows from \eqref{ylv} that the multiplicities of $0$ and $-1$ of $\E((x-\alpha)f)$ are 
the same as those of  $g$. Hence the (simple) zeros of $g$ separate the zeros of 
$\E((x-\alpha)f)$ except possibly at $0,-1$, and we conclude that $\E((x-\alpha)f)$ has only 
simple zeros.
\end{proof}
\begin{lemma}\label{e0e1}
For all integers $n \geq 1$ we have 
$$
(x+1)\E(x^n) = x\E((x+1)^n).
$$
\end{lemma}
\begin{proof}
We may write 
$$
x^n = \sum_{k=1}^na_k \binom x k, 
$$
where $a_k \in \R$. Thus 
\begin{eqnarray*} 
\E((x+1)^n) &=& \sum_{k=1}^n a_k \E [ \binom x k + \binom x {k-1} ] \\
           &=&   \sum_{k=1}^n a_k (x^k + x^{k-1}) \\
           &=& (x+1)x^{-1}\E(x^n). 
\end{eqnarray*} 
\end{proof}

Let $f$ and $g$ be standard real-rooted polynomials of degree $n$ and   
let the zeros of  $f$ and $g$ be $\alpha_1 \leq \alpha_2 \leq \cdots \leq \alpha_n$ and 
$\beta_1 \leq \beta_2 \leq \cdots \leq \beta_n$ respectively. We write $f \leq g$, if 
$\alpha_i \leq \beta_i$ for $1 \leq i \leq n$.

\begin{theorem}\label{strongalternate} 
Suppose that $f$ and $g$ are $[-1,0]$-rooted with $f \leq g$. Then 
$\E(f)$ and $\E(g)$ are $[-1,0]$- and simple-rooted,  with $\E(f) \ll \E(g)$.
\end{theorem}

\begin{proof}
By Lemma \ref{rr} and induction we only have to show that $\E(f) \ll \E(g)$. 
If $f$ and $g$ have the same zeros except for one, i.e., 
$f = (x-\alpha)h$ and $g=(x-\beta)h$, where $\alpha < \beta$, then 
$$
\E(g)= \E(f) - (\beta-\alpha)\E(h),
$$
and since $\E(h)$ interlaces $\E(f)$ we have $\E(f) \ll \E(g)$ by 
Lemma \ref{positivesum}. 

Now, suppose that $f$ and $g$ are $[-1,0]$-rooted polynomials of degree $n$ such 
that $f \leq g$. Then there are $[-1,0]$-rooted polynomials $\{h_i\}_{i=0}^M$ with 
$$
(x+1)^n=h_0 \leq h_1 \leq \cdots \leq h_M=x^n,
$$
such that $f,g \in \{h_i\}_{i=0}^M$ and $h_{i-1}$ and $h_i$ only differ in one zero 
for $1 \leq i \leq n$. 
We therefore have 
$$
\E(h_0) \ll \E(h_1) \ll \cdots \ll \E(h_M),
$$ 
and since $\E(h_0) \ll \E(h_M)$, by Lemma \ref{e0e1}, the theorem follows from  
Lemma \ref{weaktransitivity}.  
\end{proof}
A consequence of Theorem \ref{strongalternate} is that if $\{f_i\}_{i=1}^m$ is a sequence 
of standard 
$[-1,0]$-rooted polynomials of the same degree $d$, then by Lemma \ref{positivesum} and 
Theorem \ref{strongalternate}, the image under $\E$ of any non-negative sum 
$F=\sum_{i=1}^m \mu_i f_i$ will be $[-1,0]$-rooted with 
$$
\E((x+1)^d) \ll \E(F) \ll \E(x^d).
$$
It is easy to see that a standard polynomial $f$ of degree $d$ is $[-1,0]$-rooted if 
and only if $f$ can be written as 
$$
f(x) = (x+1)^dg(\frac x {x+1}), 
$$
where $g$ is a standard and $(-\infty,0)$-rooted. On the other hand, since 
$x^i(x+1)^{d-i}$ is $[-1,0]$-rooted we have that $F$ can be written as a non-negative sum 
of standard $[-1,0]$-rooted polynomials of degree $d$ if and only if   
$$ 
F(x) = \sum_{i=0}^{d} a_ix^i(x+1)^{d-i},
$$

where $a_i \geq 0$. This proves Theorem \ref{strangepf}.

\section{$t$-stack sortable permutations}\label{tstack}
For relevant definitions regarding $t$-stack sortable permutations we refer the 
reader to \cite{bonastack}.
Let $W_t(n,k)$ be the number of $t$-stack sortable permutations in the 
symmetric group, $\sym_n$, with $k$ descents, and let 
$$
W_{n,k}(x) = \sum_{k=0}^{n-1} W_t(n,k)x^k.
$$
Recently, B{\'o}na \cite{bonacor,bonasym} showed that 
for fixed $n$ and $t$ the numbers $\{W_t(n,k)\}_{k=0}^{n-1}$ form a unimodal sequence.  
When $t=n-1$ and $t=1$ we get the 
Eulerian and the Narayana numbers (see \cite{sulanke} and \cite[Exercise 6.36]{stanley2}), 
respectively. These 
are known to be $PF$-sequences and B\'{o}na \cite{bonastack,bonasym} has raised the question 
if this is true for general $t$. Here we will settle the problem to the affirmative for $t=2$ 
and $t=n-2$.

The numbers $W_2(n,k)$ are surprisingly hard to determine despite  
their compact and simple form. It was recently shown that 
$$
W_2(n,k)= \frac {(n+k)!(2n-k-1)!}{(k+1)!(n-k)!(2k+1)!(2n-2k-1)!}.
$$
See \cite{bousquet,gire,gouldenwest,jacquard} for proofs and more 
information on $2$-stack sortable permutations. 

From the case $r=0$ in Lemma \ref{factor} and the identity
$$
\sum_{k=0}^n\binom{2n-k-1}{n-1}\binom n k x^k= (-1)^n\sum_{k=0}^n\binom{-n}{k}
\binom n k (-x)^{n-k}, 
$$ 
it follows that $\binom{2n-k-1}{n-1}$ is an $n$-sequence.
\begin{theorem}
For all $n \geq 0$ the sequence $\{W_2(n,k)\}_{k=0}^{n-1}$, which records $2$-stack 
sortable permutations by descents, is $PF$.
\end{theorem}
\begin{proof}
We may write $W_2(n,k)$ as 
$$
W_2(n,k) = \frac{\binom {2n-k-1}{n-1} \binom{n+k}{n-1} \binom{2n}{2k+1}}
                {n^2\binom{2n}{n}}.
$$

A simple consequence of the notion of $PF$-sequences reads as follows: If 
$\{a_i\}_{i\geq 0}$ is $PF$ then 
so is $\{a_{ki}\}_{i\geq 0}$, where $k$ is any positive integer. 
Applying this to 
the polynomial $x(1+x)^{2n}$ we see that $\sum_k\binom{2n}{2k+1}x^k$ is 
real-rooted. Therefore the polynomial, 
$$
\sum_{k=0}^{n-1}\binom{n+k}{n-1} \binom{2n}{2k+1}x^k =
\sum_{k=0}^{n-1}\binom {2n-k-1}{n-1}\binom{2n}{2k+1}x^{n-1-k},
$$
is real-rooted. Another application 
of Lemma \ref{factor} gives that $W_{n,2}(x)$ is real-rooted. 
\end{proof}

It is easy to see that a permutation $\pi \in \sym_n$ is 
$(n-2)$-stack sortable if and only if it is not of the 
form $\sigma n1$. Thus the generating function satisfies 
$$
xW_{n,n-2}(x)= A_n(x)-xA_{n-2}(x),
$$
where $A_n(x)$ is the $n$th Eulerian polynomial.     
\begin{theorem}\label{kyrv}
For all real  numbers $t > -2$ and integers $n > 2$, the polynomial 
$$
A_n(t,x) =  A_n(x)+txA_{n-2}(x),
$$
is real- and simple-rooted. Moreover, $A_n(t,x)/x$ strictly interlaces 
$A_{n+1}(t,x)/x$ for $-2<t \leq 3$. 
\end{theorem} 
\begin{corollary}
For all $n \geq 2$ we have that $\{W_{n-2}(n,k)\}_{k=0}^{n-1}$ is $PF$.  
Moreover, $W_{n,n-2}(x)$ strictly interlaces $W_{n+1,n-1}(x)$. 
\end{corollary}
\begin{proof}[Proof of Theorem \ref{kyrv}]
It is well known that 
$A_{n-1}(x) \ll xA_{n-2}(x)$ and $A_{n-1}(x) \preceq A_n(x)$. So by 
Lemma \ref{positivesum} we have that $A_n(t,x)$ is real- and simple-rooted for 
$t \geq 0$. However, when $t <0$ a similar argument does not apply.

Let $E_n(t,x) = A_n(t, \frac x  {1+x})$. Then 
$$
E_n(t,x)=E_n(x)+tx(1+x)E_{n-2}(x),
$$
where the coefficient to $x^k$ in $E_n(x)$ counts the number of 
surjections $\sigma :  [n] \rightarrow [k]$, see \cite{brentithesis,wagnernc}. 
These polynomials satisfy the recursion:
$$
E_n(x)=x \frac{d}{dx}((1+x)E_{n-1}(x)),
$$
with initial condition $E_1(x)=x$. Thus, if we let 
$G_n(x)=E_{n+1}(x)/x$ we have the following recursion:
\begin{equation}\label{gen}
G_n(x)=\frac{d}{dx}(x(1+x)G_{n-1}(x)),
\end{equation}
with $G_0(x)=1$. Obviously $G_n(x)$ is real- and simple-rooted. If 
we apply \eqref{gen} two times we get the equation:
\begin{multline*}
G_n(x)= (1+6x+6x^2)G_{n-2}(x)+3x(1+2x)(1+x)G'_{n-2}(x)+ \\ 
         x^2(1+x)^2G''_{n-2}(x),
\end{multline*}
and for $G_n(t,x):= G_n(x)+tx(1+x)G_{n-2}(x)$ we have 
\begin{multline*}
G_n(t,x) =  (1+(6+t)x+(6+t)x^2)G_{n-2}(x)+3x(1+2x)(1+x)G'_{n-2}(x)+\\ 
          x^2(1+x)^2G''_{n-2}(x).
\end{multline*}
To apply Theorem \ref{hyper} we need     
show that for all $\xi \in \R$ and $-2<t<0$ the polynomial 
$$
F(\xi,z):=(1+(6+t)\xi+(6+t)\xi^2)+3\xi(1+2\xi)(1+\xi)z 
         + \xi^2(1+\xi)^2z^2
$$
is real-rooted. The discriminant of $F(\xi,z)$, 
$$
\Delta(F(\xi,z))=\xi^2(1+\xi)^2(2+t+(3-t)(1+2\xi)^2),
$$
is non-negative when $-2 \leq  t \leq 3$, so $F(\xi,z)$ real-rooted for these $t$. 
Since all the  
$Q_k$s are standard it is easy to see that condition (III) in the statement of 
Theorem \ref{hyper} is satisfied. 
Moreover, $1+(6+t)x+(6+t)x^2$ strictly interlaces $3x(1+2x)(1+x)$ when 
$t > -2$ so Theorem 
\ref{hyper} applies. Since $G_n$ strictly interlaces $G_{n+1}$ we have by 
Theorem \ref{hyper} and Corollary \ref{obreschkor} that 
$\phi_F(G_n)$ strictly interlaces $\phi_F(G_{n+1})$. Thus 
$A_n(t,x)$ strictly interlaces $A_{n+1}(t,x)$. 
\end{proof} 

\section{$q$-Eulerian and $W$-Eulerian polynomials}\label{qeuler}
A $q$-analog of the Eulerian polynomials was introduced and studied in \cite{foata} and 
further studied in \cite{brentiplethysm}. It is defined by 
$$
A_n(x;q):= \sum_{\pi \in \sym_n}x^{\exc(\pi)}q^{c(\pi)},
$$
where $\c(\pi)$ and $\exc(\pi)$ denotes the number of {\em cycles} and {\em excedances} 
in $\pi$ respectively. These 
polynomials satisfy the recursion 
$$
A_{n+1}(x;q) = (nx+q)A_n(x;q)-x(x-1)\frac{\partial}{\partial x} A_n(x;q),
$$
with initial condition $A_0(x;q):=1$. See \cite{brentiplethysm} for a proof. 
The following theorem appears in \cite{brentiplethysm}. 
\begin{theorem}\label{brentitheorem}
Let $q \in \R$, $q>0$. Then the polynomials $A_n(x,q)$ have only real 
non-positive simple zeros. 
\end{theorem}
Brenti also makes the following conjecture:
\begin{conjecture}\label{brenticonjecture}
Let $n,m \in \N$. Then $A_n(x;-m)$ has only real zeros.
\end{conjecture}
In what follows we will prove this conjecture using multiplier $n$-sequences. For 
$n \in \N$ define the polynomials $E_n(x;q)$ by:
$$
E_n(x;q):= (1+x)^nA_n(\frac x {1+x};q).
$$
It is clear that $E_n(x;q)$ is real-rooted if and only if $A_n(x;q)$ is 
real-rooted. These polynomials satisfy a somewhat easier recursion. Namely, 
\begin{equation}\label{E-recursion}
E_{n+1}(x;q) = (1+x)\{ qE_n(x;q)+x\frac{\partial}{\partial x}E_n(x;q) \},
\end{equation}
with initial condition $E_0(x;q)=1$. Now, for $q \in \R$ let 
$\Gamma_q : \R[x] \rightarrow \R[x]$ be the linear operator defined 
by $\Gamma_q(f(x)) = qf(x) + xf'(x)$. Since $\Gamma_q(x^n)= (q+n)x^n$ we 
may apply Corollary \ref{gammaq}. 
\begin{theorem}
Let $q \in \R$ and $n \in \N$. If $q \geq 0$, $n \leq -q$ or $q \in \Z$ then 
$E_n(x;q)$ has only real zeros. 
\end{theorem}

\begin{proof} 
We may write \eqref{E-recursion} as 
$$
E_{n+1}(x;q)=(x+1)\Gamma_q[E_n(x;q)].
$$
The cases $q \geq 0$ and $n \leq -q$ follow from Corollary \ref{gammaq} by 
induction. We may therefore assume that $q=-m$ is a negative integer. 
We claim that $\deg E_n(x;q) = n$ if $n \leq m$ and $\deg E_n(x;q)= m$ if 
$n \geq m$. From this the real-rootedness follows by Corollary \ref{gammaq} 
and induction. The case $n \leq m$ is clear since 
$\Gamma_q[x^{n-1}]=-(m-n+1)<0$. The case $n>m$ also follows by induction. Suppose 
that $n \geq m$ and that $\deg E_n(x;q)=m$. Then by the recursion we 
have that $\deg E_{n+1}(x;q) \leq m+1$. Moreover, since 
$\Gamma_q[x^m] = 0$ we have that  $\deg E_{n+1}(x;q) \leq m$. Let 
$a \neq 0$ be the coefficient to $x^m$ of $E_{n}(x;q)$. Then the 
coefficient to $x^m$ of $E_{n+1}(x;q)$ is $a \Gamma_q[x^{m-1}] = -a$, 
so $\deg E_{n+1}(x;q) = m$, and the thesis follows.
\end{proof}

The Eulerian polynomial, $P(W,x)$, of a finite Coxeter group $W$ is 
the polynomial, 
$$
P(W,x)= \sum_{\sigma \in W}x^{d_W(\sigma)},
$$
where $d_W(\sigma)$ is the number of $W$-descents of $\sigma$, see \cite{brentieuler}. 
This polynomial is also the generating function for the $h$-vector of the Coxeter complex 
associated to $(W,S)$. For Coxeter groups of type $A_n$ we have that 
$P(A_n,x)=A_n(x)/x$, the shifted Eulerian polynomial. Also, for Coxeter groups 
of type $B_n$ it is known, see \cite{brentieuler}, that $P(B_n,x)$, has only real zeros. 
It is easy to see that $P(W_1 \times W_2,x) = P(W_1,x)P(W_2,x)$ for finite Coxeter groups 
$W_1$ and $W_2$. Also, the real-rootedness can be checked ad hoc for 
the exceptional groups. Thus, by the classification of finite irreducible Coxeter 
groups, to prove that $P(W,x)$ has only real zeros for all finite Coxeter 
groups it suffices to prove that $P(D_n,x)$ is real-rooted for Coxeter groups 
of type $D_n$. The real-rootedness of $P(D_n,x)$ is conjectured by Brenti in 
\cite{brentieuler}. It is known that the Eulerian polynomials of type $A_n,B_n$ and 
$D_n$ are related by, see \cite{brentieuler,reiner,stembridge}: 
$$
P(D_n,x)= P(B_n,x)-n2^{n-1}xP(A_{n-1},x). 
$$ 
This relationship was first noticed by Stembridge \cite{stembridge}.  
One step towards proving the real-rootedness of $P(D_n,x)$ is to learn 
more about the relationships between the zeros of $P(B_n,x)$ and $P(A_n,x)$.        
    
Brenti \cite{brentieuler} introduced a $q$-analog of $P(B_n,x)$  
\begin{equation}\label{B-q}
B_n(x;q) = \sum_{\sigma \in B_n} q^{N(\sigma)}x^{d_B(\sigma)},
\end{equation}
where $d_B(\sigma)$ is the number of $B_n$-descents of $\sigma$ and $N(\sigma)$ is the 
number of negative entries of $\sigma$, see \cite{brentieuler}. He proved that 
\begin{equation}\label{B_n}
\sum_{i \geq 0}( (1+q)i+1)^nx^i = \frac{B_n(x;q)}{(1-x)^{n+1}},
\end{equation}
and that $B_n(x;q)$ is real- and simple-rooted for all $q \geq 0$. 
Suppose that $f(i)$ is a polynomial in $i$ of degree $d$, then the 
polynomial 
$W(f)$ is defined by 
$$
\sum_{i \geq 0}f(i)x^i = \frac{W(f)(x)}{(1-x)^{d+1}},
$$
One can show, see \cite{brentithesis}, that $\E(f)$ and $W(f)$ are related by:
\begin{equation}\label{label}
\E(f)(x) = (1+x)^{\deg(f)} W(f)(\frac x {1+x}). 
\end{equation}
It follows that $W(f)$ has only real non-positive roots if and only if $\E(f)$ is 
$[-1,0]$-rooted. 

Since $((1+q)i+1)^n$ is a $[-1,0]$-rooted polynomial in $i$ for any $q \geq 0$ it follows 
from e.g. Theorem \ref{strangepf} that 
$B_n(x;q)$ is real-rooted in $x$ for any fixed $q \geq 0$. It is 
natural to generalize $B_n(x;q)$ to have $n+1$ parameters as 
$B_n(x;{\bf q}) := W\big(\prod_{i=0}^n \big( (1+q_i)x +1 \big)\big)$. This 
polynomial has a nice combinatorial interpretation:
\begin{theorem}\label{multieuler}
For all $n \in \N$ we have:
$$
B_n(x,{\bf q}) = \sum_{\sigma \in B_n} 
q_1^{\chi_1(\sigma)}q_2^{\chi_2(\sigma)} \cdots q_n^{\chi_n(\sigma)}t^{d_B(\sigma)},
$$
where 
$$
\chi_i(\sigma) = \begin{cases}
                   1  \ \ \mbox{ if } \sigma_i < 0, \\
                   0  \ \ \mbox{ if }  \sigma_i > 0. 
                   \end{cases}
$$
\end{theorem} 
\begin{proof}
The proof is an obvious generalization of the proof of Theorem 3.4 of \cite{brentieuler}. 
\end{proof}
Note that this theorem gives a semi-combinatorial interpretation of the 
$W$-transform of any $[-1,0)$-rooted polynomial. 
\begin{corollary}
Let $n \in \N$ and let $q_1,q_2, \ldots, q_n$ be  non-negative real numbers. Then 
$B_n(x;{\bf q})$ has only real and simple zeros.
\end{corollary}

We need the following lemma on the degree of $W(f)$. 
\begin{lemma}\label{degw}
Let $f \in \R[x]$. Then 
$$
\deg W(f) = \deg f -  \mult (-1, \E(f)).
$$
Moreover, $\mult (-1, \E(f))$ is equal to the maximal integer $k$ such that 
$(x+1)(x+2)\cdots (x+k)$ divides $f$. 
\end{lemma}
\begin{proof}
Since $\deg \E(f) = \deg f$ for all $f$ we have by \eqref{label} that 
$\deg W(f) = \deg f -  \mult (-1, \E(f))$. If we expand $f$ in the basis 
$\{ \binom {-x-1} i\}$ as: 
\begin{eqnarray*}
f(x) &=& \sum_{i\geq 0}(-1)^ia_i \binom {-x-1} i,\\
     &=& \sum_{i\geq 0}\frac {a_i}{i!} (x+1) \cdots (x+i),
\end{eqnarray*}
we have by Lemma \ref{symE} that 
$$
\E(f)(x) = \sum_{i\geq 0}a_i(x+1)^i,
$$
and the lemma follows.
\end{proof}
We now have more precise knowledge of the location of the zeros of $B_n(x;q)$ for 
any given $q \geq 0$.   
\begin{theorem}
Let $0 <q <t \in \R$ and $n >0$ be an integer. Then  
$$
B_n(x;0) \preceq B_n(x;t) \ll B_n(x;q) \ll xB_n(x;0),
$$
where the three first polynomials have no common zeros.
\end{theorem}
\begin{proof}
Let $0 < r < s < 1$. Then by the proof of Lemma \ref{rr} we have 
\begin{multline*} 
\E(x^n) \ll \E(x(x+r)^{n-1}) \ll_{strict} \E((x+r)^n) 
\ll_{strict} \E((x+r)^{n-1}(x+s)) \ll \\ \E((x+s)^n) \ll_{strict} 
\E((x+s)^{n-1}(x+1)) \ll \E((x+1)^n),
\end{multline*}
where $\ll_{strict}$ means strictly alternating left of. Since 
$(x+1)\E(x^n) = x\E((x+1)^n)$ this implies 
$$
\E(x^n) \ll_{strict} \E((x+r)^n) \ll_{strict} \E((x+s)^n) \ll_{strict} \E((x+1)^n).
$$
Now since 
$$
B_n(x;q)= (q+1)^nW((x+\frac{1}{1+q})^n) = 
(q+1)^n(1-x)^n\E( (x+\frac{1}{1+q})^n)(\frac{x}{1-x}),
$$
we see by Lemma \ref{degw} that $\deg B_n(x;0) = n-1$ and $\deg B_n(x;q) = n$ if $q \neq 0$. 
Moreover, the alternating property is preserved under the operation \eqref{label} 
and the theorem follows.
\end{proof}
It follows from \eqref{B-q} that $P(B_n,x)= B_n(x;1)$ and $P(A_n,x)=B_n(x;0)$.
\begin{corollary}
For all integers $n \geq 1$ we have that 
$P(A_n,x)$ strictly interlaces $P(B_n,x)$.
\end{corollary} 
Since $P(A_n,x) \ll xP(A_{n-1},x)$ and $P(A_n,x) \preceq P(B_n,x)$, we have by Lemma 
\ref{positivesum} that for all $t \geq 0$ the polynomial $P(B_n,x) + txP(A_{n-1},x)$ 
is real-rooted. Unfortunately a similar argument does not apply when $t<0$.   

One can extract more from \eqref{B_n}. 
Brenti \cite{brentieuler} proved that the polynomial 
$$
\sum_{\sigma \in B_n, N(\sigma) \in \{k,n-k\}}x^{d_B(\sigma)},
$$
is real-rooted for all choices of $0 \leq k \leq n$. 
Using Theorem \ref{strongalternate} we can extend this result to:
\begin{corollary}\label{negativeset}
Let $S$ be any subset of $[0,n]$. Then the polynomial 
$$
P(B_n,S;x):= \sum_{\sigma \in B_n, N(\sigma) \in S} x^{d_B(\sigma)},
$$
has only real and simple zeros.
\end{corollary}
\begin{proof}
Comparing the coefficient of $q^i$ in both sides of \eqref{B_n} we see that 
$P(B_n,S;x)= W(f_n(S;x))$ where 
$$
f_n(S;x)=\sum_{s \in S} \binom n s x^s(x+1)^{n-s}.
$$
So the theorem follows from Theorem \ref{strangepf}.
\end{proof}
One instance of Theorem \ref{negativeset} is particularly interesting. Recall 
that a Coxeter group of type $D_n$ is isomorphic to the subgroup
$$
D_n= \{ \sigma \in B_n : \ \ 2 \mid N(\sigma)\}
$$
Hence, we have the following corollary
\begin{corollary}\label{wrongdescent}
For all $n \in \N$ the polynomial 
$$
\sum_{\sigma \in D_n} x^{d_B(\sigma)}
$$
has only real and simple zeros. 
\end{corollary}
Note that the above polynomial is not $P(D_n,x)$, since $B_n$-descents and 
$D_n$-descents are not the same.

\section{The $h$-vector of a family of simplicial complexes 
defined by Fomin and Zelevinsky}\label{hvector}

Fomin and Zelevinsky \cite{fomin} recently associated to any finite Weyl group 
$W$ a simplicial complex $\Delta_{FZ}(W)$. For the classical Weyl groups 
these polynomials are given by 
\begin{eqnarray*}
h(\Delta_{FZ}(A_{n-1}),x)&=& \frac {1} {n}\sum_{k=0}^{n-1}\binom {n} {k} \binom {n} {k+1}x^k,\\
h(\Delta_{FZ}(B_n),x)&=& \sum_{k=0}^n \binom n k \binom n kx^k,\\
h(\Delta_{FZ}(D_n),x)&=& h(\Delta_{FZ}(B_n),x) - nx h(\Delta_{FZ}(A_{n-2}),x).
\end{eqnarray*}
It is known that the $h$-polynomials corresponding to  $A_n$ and $B_n$ have only real zeros. 
We will here show that so has $h(\Delta_{FZ}(D_n),x)$. 
\begin{theorem}\label{mainweyl}
Let $\alpha, \beta \in \R$ be such that $\alpha \geq 0, 2\alpha + \beta>0$ and let 
$n \geq 2$ be an integer. 
Then the polynomial 
$$
F_n(\alpha,\beta):=\alpha h(\Delta_{FZ}(B_n),x) + \beta nx h(\Delta_{FZ}(A_{n-2}),x),
$$
is real- and simple-rooted. Moreover, $h(\Delta_{FZ}(B_{n-1}),x)$ strictly interlaces 
$F_n(\alpha,\beta)$ if $\alpha>0$ and strictly alternates left of $F_n(\alpha,\beta)$ 
if $\alpha=0$. 
\end{theorem}
\begin{corollary}
Let $W$ be a finite Weyl group. Then $h(\Delta_{FZ}(W),x)$ has only real and 
simple zeros. 
\end{corollary}
\begin{proof}
For the exceptional Weyl group one can check the real-rootedness ad hoc, see 
\cite{reinerwelker}. The other cases follows from 
Theorem \ref{mainweyl}.
\end{proof}
The {\em Hadamard product} of two polynomials 
\begin{eqnarray*}
p(x)&=& a_0 + a_1x + \cdots + a_mx^m \\
q(x)&=& b_0 + b_1x + \cdots + b_nx^n
\end{eqnarray*}
is the polynomial 
$$
(p \star q)(x) = a_0b_0 + a_1b_1x + \cdots + a_Nb_Nx^N, 
$$
where $N=\min(m,n)$. M\'alo proved that if the zeros of $p$ are real and the zeros 
of $q$ are real and of the same sign then the zeros of $p \star q$ are real as well. 
This also follows from Theorem \ref{schurmult} since 
$p \star q = \Gamma [p S q]$ where 
$\Gamma$ is the multiplier sequence $\{ \frac{1}{k!} \}_{k=0}^{\infty}$. It is 
known, see e.g. \cite{garloffwagner}, that if $f$ has only real zeros then 
all zeros of $\Gamma[f]$ are real and simple except for possibly at the origin.   
\begin{proof}[Proof of Theorem \ref{mainweyl}]
We may write $F_n(\alpha,\beta)$ as 
$$
F_n(\alpha,\beta)=\alpha (x+1)f +(2\alpha + \beta)g,
$$
where $f= (x+1)^{n-1} \star (x+1)^{n-1}$ and $g=\big(x(x+1)^{n-1}\big) \star (x+1)^{n-1}$. 

By the discussion before this proof we have that for all real choices 
of $\gamma, \delta \in \R$ the polynomial 
$$
\gamma f + \delta g = \big( (\gamma+ \delta x)(x+1)^{n-1}\big)\star (x+1)^{n-1},
$$
is real- and simple-rooted. By the Obreschkoff theorem we infer that 
$f$ strictly alternates left of $g$. Now, since $f \preceq (x+1)f$ and 
$f \ll g$ we know by Lemma \ref{positivesum} that $f$ either interlaces 
or alternates left of  $F_n(\alpha,\beta)$ for all $\alpha,\beta \in \R$ such that 
$\sgn (\alpha) = \sgn (2\alpha+\beta)$. Moreover, since $g$ and $f$ have no common 
zeros nor does $F_n(\alpha,\beta)$ and $f$ (provided that $2\alpha + \beta \neq 0$).   
\end{proof}

\bigskip 

\noindent {\bf Acknowledgments.}  The author would like to thank Francesco Brenti 
for interesting discussions and helpful suggestions on the topic of this paper during 
the author's stay at Universit\'a di Roma ``Tor Vergata''. Thanks to an anonymous 
referee for valuable comments on the presentation of this paper.

\bibliography{realpol}

\end{document}